\documentclass[10pt,a4paper]{amsart}
\usepackage{amssymb}

\usepackage{color}
\usepackage{hyperref}
\usepackage{dsfont}

\usepackage{tikz}

\newtheorem{thm}{Theorem}

\newtheorem{cor}{Corollary}

\newenvironment{pf}{{\it Proof:}\quad}{\hfill$QED$}
\newcommand{\R}{\mathbb R}

\newcommand{\C}{\mathbb C}
\newcommand{\Z}{\mathbb Z}

\begin{document}

\title{Borel class and Cartan involution}
 \author{Thilo Kuessner}
   \email{thilo.kuessner@math.uni-augsburg.de}

\maketitle

\begin{abstract}Let $\rho\colon \pi_1M^3\to PGL(n,\C)$ be a representation of a $3$-manifold group to the projective linear group. In this note we prove that the Borel class of $\rho$ is preserved up to sign under Cartan involution.

For representations to $PGL(3,\C)$ this is implied by the
more general result of E.\ Falbel and
Q.\ Wang from \cite{fw}, however our proof appears to be much shorter for that special case. 
\end{abstract}

\section{Introduction}Representations of fundamental groups of (ideally triangulated) $3$-manifolds to $PGL(3,\C)$ are determined by so-called decorations, i.e.,  equivariant maps from the 0-skeleton of the universal covering to $Flag(\C P^2)$. Bergeron, Falbel and Guilloux introduced in \cite{bfg} invariants of $4$-tuples of flags, which for a given decoration determine an element $\beta(T)$ in the pre-Bloch group ${\mathcal P}(\C)$ for each tetrahedron $T$. Summing up the $\beta(T)$ over all simplices in the triangulation defines an element in the Bloch group, whose Bloch-Wigner dilogarithm divided by $4$ is said to be the volume of the representation.

In \cite{fw}, Falbel and Wang considered the effect of the involution $A\to (A^T)^{-1}$ (or on the level of decorations of going from flags to their duals) on these invariants. Using ingenious computations in the pre-Bloch group they expressed in \cite[Proposition 3.6]{fw} the
difference of the values of $\beta(T)$, for a decoration and its dual,
as the sum of the Fock-Goncharov triple ratios associated to the induced decorations of the four faces of the simplex $T$. When summing up over all 3-simplices, faces occuring twice with opposite orientations yield canceling difference terms. In particular, when the boundary consists of cusps only, the authors obtain that the element in the Bloch group (and hence the volume) is the same for the decoration and its dual. This
proves that the involution $A\to (A^T)^{-1}$ of $PGL(3,\C)$ preserves the volume and the Cartan involution $A\to(\overline{A}^T)^{-1}$ 
multiplies it by $-1$. 

In the introduction of \cite{fw} the authors point out that this implies the same invariance properties for the Borel class and they express the belief that an analogous invariance under duality should hold for representations to $PGL(n,\C)$ with $n\ge 4$, and that also for non-toroidal boundary one should then have an exact coboundary term relating the two cocycle representatives coming from the decoration and its dual.

In this paper we prove the first part of this conjecture with a short argument which in particular also gives a much shorter proof for the case $n=3$.

What we actually prove is that the Cartan involution preserves the volume and the Chern-Simons invariant of representations of fundamental groups up to a sign depending on dimension. It follows from the van Est theorem that all invariants coming from the continuous cohomology group $H^3_c(PGL(n,\C))$ must be multiples of each other, so the invariance of the volume implies the invariance of the Borel class and also of the volume defined in \cite{bfg}.

Our proof via secondary characteristic classes does not yield information on the pre-Bloch group elements associated to the individual simplices. So it remains open whether, as conjectured in \cite{fw},
also for $n\ge 4$ a simplex-wise formula analogous to \cite[Proposition 3.6]{fw} relates the two values of $\beta(T)$ for a decoration and its dual.

In \hyperref[pgl3cbfg]{Section \ref*{pgl3cbfg}} we include a proof that the volume of \cite{bfg} agrees up to multiplication by $-\frac{1}{4}$ with the volume defined via Cheeger-Chern-Simons classes. This seems not to have been known (see \cite[Remark 9.\ 19]{gtz}), although the proof is actually a straightforward application of the methods introduced in \cite{gtz}.
\section{Preliminaries}Let $M$ be a compact, orientable, (2k-1)-dimensional manifold, possibly with boundary $\partial M$. Let $\rho\colon\pi_1M\to SL(n,\C)$ be a boundary-unipotent representation of its fundamental group, where boundary-unipotent means that the image of the fundamental group 
of $\partial M$ can be conjugated into the group $N(n,\C)\subset SL(n,\C)$ of upper triangular matrices with 1's on the diagonal. 

Let $SL(\C)=\bigcup_nSL(n,\C)$ and $N(\C)=\bigcup_n N(n,\C)$ be the increasing unions with respect to the canonical embeddings as block matrices. It is well known that the image of the relative fundamental class in group homology
$$(B\rho)_*\left[M,\partial M\right]\in H_{2k-1}(BSL(\C)^\delta,BN(\C)^\delta;Z)$$ has a preimage in $H_{2k-1}(BSL(\C)^\delta;\Z)$, which we will denote by the same name.

Let $\hat{c}_k\in H^{2k-1}(BSL(\C)^\delta,\C/\Z)$ be the universal Cheeger-Chern-Simons class. Then $$\langle \hat{c}_k,(B\rho)_*\left[M,\partial M\right]\rangle\in\C/\Z$$
is the Cheeger-Chern-Simons invariant of the flat bundle with monodromy $\rho$. It has become common, at least in the context of $3$-manifolds, to denote its real and imaginary part as the Chern-Simons invariant and the volume of the representation, respectively. See \cite[Section 2]{gtz}.

We remark that this approach implies that volume and Chern-Simons invariant are constant on components of the 
character variety of boundary-unipotent representations, because they can be computed from the fundamental class 
in the discrete subset $H_{2k-1}(BSL(\C)^\delta;\Z)\subset H_{2k-1}(BSL(\C)^\delta;\R)$ which can not vary continuously. This local constantness is hard to see from other definitions, compare \cite{bas} and \cite{lod} for other proofs in the case of $\partial M=\emptyset$.

\section{Invariance of Cheeger-Chern-Simons-classes}
\begin{thm}\label{invpol}
Let $M$ be a compact, orientable $(2k-1)$-manifold. Let $E_\rho\to M$ be a trivializable, flat $SL(n,\C)$-bundle with monodromy $\rho\colon\pi_1M\to SL(n,\C)$,  let $\rho^\vee\colon\pi_1M\to SL(n,\C)$ be the dual of $\rho$ under Cartan involution.
Then $$vol(\rho^\vee)=(-1)^{k+1}vol(\rho),\ \  CS(\rho^\vee)=(-1)^kCS(\rho)$$
\end{thm}
\begin{pf}
Recall (e.g., from \cite{bas}) that $\hat{c}^k(E_\rho)$ can be computed as 
$$\hat{c}^k(E_\rho)=\frac{1}{(2\pi)^k}\int _Ms^*(A_kP_k(\theta))\in \C/\Z$$
for the flat connection $\theta\in\omega^1(E_\rho,\mathfrak{sl}(n,\C))$, with the constant $A_k=(-1)^{k-1}\frac{k!(k-1)!}{2^{k-1}(2k-1)!}$ and the invariant polynomial $P_k(A)=Tr(A^k)$, and with a section $s\colon M\to E_\rho$ which exists by trivializability of the bundle. (The $mod\ \Z$-indeterminacy comes from the possible choice of different trivializations.)

For any group homomorphism $\phi\colon SL(n,\C)\to SL(n,\C)$ we obtain a trivializable, flat bundle with the "same" section $s\colon M\to E_{\phi\circ\rho}$, the 
monodromy $\phi\circ\rho$ and flat connection $\phi_*\theta$, where $\phi_*=D_{\mathds 1}\phi$ denotes the induced homomorphism of Lie algebras. Thus we have
$$\hat{c}^k(E_{\phi\circ\rho})=\frac{1}{(2\pi)^k}\int _Ms^*(A_kP_k(\phi_*\theta)).$$

The Cartan involution is given by $\phi(A)=(\overline{A}^T)^{-1}$ and its differential at ${\mathds 1}$ is $$\phi_*(X)=-\overline{X}^T.$$
Clearly
$$Tr((-\overline{X}^T)^k)=(-1)^k\overline{Tr(X^k)},$$
which implies $$P_k(\phi_*\theta)=(-1)^k\overline{P_k(\theta)},$$
from which the claim follows because volume and Chern-Simons invariant are defined by the imaginary and real part of $\hat{c}_2$, respectively.
\end{pf}

\begin{cor}\label{cscor}For a compact, orientable $3$-manifold $M$ and a representation $\rho\colon \pi_1M\to SL(n,\C)$ we have $vol(\rho^\vee)=-vol(\rho),\  CS(\rho^\vee)=CS(\rho)$.\end{cor}
\begin{pf}$SL(n,\C)$-bundles over $3$-manifolds are trivializable, thus \hyperref[invpol]{Theorem \ref*{invpol}} applies.\end{pf}

\section{The Borel class}\label{borelclass}
The van Est theorem implies $H^3_c(PGL(n,\C);\R)\cong H^3_{dR}(SU(n))\cong \R$. The generator of the image of $H^3_c(PGL(n,\C);\Z)$ is called the Borel class $b_3$. For representations $\rho$ of $3$-manifold groups, their Borel invariant is defined as $\langle \rho^*b_3,\left[M\right]\rangle$. See \cite{bbi} for an explicit description of the Borel invariant.

The imaginary part of the Cheeger-Chern-Simons class is a continuous $\R$-valued cohomology class and thus a nontrivial multiple of the Borel class. Hence \hyperref[cscor]{Corollary \ref*{cscor}} implies that the Borel class is preserved under the involution $A\to(A^T)^{-1}$, which has been conjectured in the introduction of \cite{fw}.

\section{The PGL(3,C) case}\label{pgl3cbfg}

One of the results in \cite{fw} is the equality $vol(\rho^\vee)=-vol(\rho)$ for the volume of representations defined in \cite{bfg}. We claim that the the volume from \cite{bfg} 
agrees with $-\frac{1}{4}$ times 
the imaginary part of the Cheeger-Chern-Simons class and thus this result is a special case of \hyperref[invpol]{Theorem \ref*{invpol}}. (Using the same reasoning 
as in \hyperref[borelclass]{Section \ref*{borelclass}} it would actually be enough to apply the van Est theorem to recover that
result from \hyperref[invpol]{Theorem \ref*{invpol}}, but we think that the precise relation between the two volumes should be of independent interest.)

We remark that this relation of volumes seems not to appear in the literature, see \cite[Remark 9.19]{gtz}, although its proof is actually a rather straightforward application of the results in \cite{gtz}.

\begin{thm}For a boundary-parabolic representation $\rho\colon\pi_1M^3\to PGL(3,{\bf C})$ of a $3$-manifold group, its volume $vol(\rho)$
as defined in \cite{bfg} agrees with $-\frac{1}{4}$ times the imaginary part of the Cheeger-Chern-Simons invariant $\hat{c}_2$.\end{thm}
\begin{pf} According to the main result of \cite{gtz}, specialized to $G=PGL(3,{\C})$, the imaginary part of the 
Cheeger-Chern-Simons invariant can be computed as follows. Fix a 
generalized ideal triangulation of $M$, let $L$ be the corresponding generalized 
triangulation of the space obtained from the universal covering $\widetilde{M}$ 
by collapsing each boundary component to a point, and let $L_0$ be its 
$0$-skeleton. Let $N\subset G$ be the group of upper triangular matrices with $1$'s on the diagonal and let 
$\overline{\xi}\colon L_0\to G/N$ be a $\rho$-equivariant map, which can for example be defined by mapping 
one vertex in each orbit arbitrarily and then extending equivariantly. (Stabilizers of points in $L_0$ are fundamental 
groups of boundary components, thus boundary-parabolicity 
guarantees that extending equivariantly yields a well-defined map to $G/N$.) Then for each simplex of $K$ the images under $\xi$ of the vertices of its lifts to $L$ are a $\rho(\pi_1M)$-orbit of
$4$-tuples $(g_0N,g_1N,g_2N,g_3N)$. Then, following \cite[Definition 5.1]{gtz}, for any $4$-tuple 
$t=(t_0,t_1,t_2,t_3)$ of nonnegative integers with $t_0+t_1+t_2+t_3=3$ one considers its ptolemy coordinate
$c_t=det(\cup_{i=0}^3\left\{g_i\right\}_{t_i})$, where $\left\{g_i\right\}_{t_i}$ means the ordered set of the first $t_i$
column vectors of $g_i$ and $\cup_{i=0}^3\left\{g_i\right\}_{t_i}$ means the matrix whose ordered column set is composed by the $\left\{g_i\right\}_{t_i}$. (This is well-defined for each simplex of $K$, i.e., does not depend on the lift to $L$.)
With this ptolemy coordinates one computes, following \cite[Section 5.3]{gtz} an element $$\hat{c}(\rho)=$$
$$(\log(c_{2001})+\log(c_{1110})-\log(c_{2010})-\log(c_{1101}),\log(c_{2100})+\log(c_{1011})-\log(c_{2010})-\log(c_{1101})$$
$$+(\log(c_{1101})+\log(c_{0210})-\log(c_{1110})-\log(c_{0201}),\log(c_{1200})+\log(c_{0111})-\log(c_{1110})-\log(c_{0201})$$
$$+(\log(c_{1011})+\log(c_{0120})-\log(c_{1020})-\log(c_{0111}),\log(c_{1110})+\log(c_{0021})-\log(c_{1020})-\log(c_{0111})$$
$$+(\log(c_{1002})+\log(c_{0111})-\log(c_{1011})-\log(c_{0102}),\log(c_{1101})+\log(c_{0012})-\log(c_{1011})-\log(c_{0102})\in\widehat{\mathcal P}(\C)$$
in the extended Bloch group $\widehat{\mathcal P}(\C)$ associated to the simplex. Summing over all simplices of the triangulation one obtains an 
element $\rho_*\left[M,\partial M\right]\in \widehat{\mathcal B}(\C)\subset \widehat{\mathcal P}(\C)$, such that application of the extended Roger's dilogarithm $R\colon \widehat{\mathcal P}(\C)\to\C/4\pi^2\Z$ to this element gives the Cheeger-Chern-Simons invariant $\hat{c}_2$. See \cite[Theorem 8.3]{gtz}.

Since we are interested only in its imaginary part we can use the relation between the extended Roger's dilogarithm and the Bloch-Wigner dilogarithm: the imaginary part of the extended Rogers' dilogarithm 
of an element in $\widehat{\mathcal P}(\C)$ is the Bloch-Wigner dilogarithm of its image in 
${\mathcal P}(\C)$.
For the above element $\hat{c}(\rho)\in \widehat{\mathcal P}(\C)$ its image in ${\mathcal P}(\C)$ is
\begin{equation}
c(\rho)= \left[\frac{c_{2001}c_{1110}}{c_{2010}c_{1101}}\right]+
\left[\frac{c_{1101}c_{0210}}{c_{1110}c_{0201}}\right]+
\left[\frac{c_{1011}c_{0120}}{c_{1020}c_{0111}}\right]+
\left[\frac{c_{1002}c_{0111}}{c_{1011}c_{0102}}\right]\in {\mathcal P}(\C).
\end{equation}
One obtains the volume associated to the simplex of $K$ by applying the Bloch-Wigner dilogarithm $D\colon {\mathcal P}(\C)\to \R$ to $c(\rho)$, and the volume of $\rho$ by summing these volumes over all simplices of the triangulation.

From $\overline{\xi}\colon L_0\to G/N$ one obtains a $\rho$-equivariant map $\xi\colon L_0\to Flag(\C P^2)$ via the obvious projection $G/N\to G/B=Flag(\C P^2)$ with $B\subset G$ the group of all upper triangular matrices. For given $\xi$, the volume of the associated flag structure is according to \cite[Section 3.6]{bfg} 
defined by summing over all simplices of the triangulation $\frac{1}{4}$ times 
the Bloch-Wigner dilogarithms of
\begin{equation}
\left[z_{12}\right]+\left[z_{21}\right]+\left[z_{34}\right]+\left[z_{43}\right]\in{\mathcal P}(\C),
\end{equation}
where the numbers in brackets are associated to a 4-tuple $(\left[x_i\right],\left[f_i\right])_{1\le i\le 4}$ of flags as follows. 

To define $z_{ij}$ define $k$ and $l$ such that the permutation $(1,2,3,4)\to(i,j,k,l)$ is even. The lines through $x_i$ form a projective line, so one has a well-defined cross ratio of four lines through $x_i$. Then $z_{ij}$ is defined as the cross ratio of the lines $ker(f_i), (x_ix_j), (x_ix_k), (x_ix_l)$.

According to \cite[Lemma 2.7]{bfg} the value of $z_{ij}$ can be computed by the formula
$$
z_{ij}=\frac{f_i(x_k)det(x_i,x_j,x_l)}{f_i(x_l)det(x_i,x_j,x_k)}.
$$

Now consider the map $G/N\to G/B=Flag(\C P^2)$. Denote $(e_1,e_2,e_3)$ the standard basis of $\C^3$ and $(x_0,f_0)=(\left[e_1\right],\left[\langle e_3,.\rangle\right])$ the standard flag in $\C P^2$. For $g_i\in G/N$ we denote $(x_i,f_i)=(g_i(x_0),g_i(f_0))$ the corresponding element in $G/B=Flag(\C P^2)$. 

An elementary linear algebra computation shows $$det(\left\{g_i\right\}_{1}\left\{g_j\right\}_{1}\left\{g_k\right\}_{1})=det (x_i,x_j,x_k)$$ and 
$$det(\left\{g_i\right\}_{2}\left\{g_j\right\}_{1})=f_i(x_j)$$
for all $i,j,k$. Thus, for example, 
$$\left[\frac{c_{2001}c_{1110}}{c_{2010}c_{1101}}\right]=\frac{f_1(x_4)det(x_1,x_2,x_3)}{f_1(x_3)det(x_1,x_2,x_4)}=\left(\frac{f_1(x_3)det(x_1,x_2,x_4)}{f_1(x_4)det(x_1,x_2,x_3)}\right)^{-1}$$
and hence $$D(\left[\frac{c_{2001}c_{1110}}{c_{2010}c_{1101}}\right])=-D(\frac{f_1(x_3)det(x_1,x_2,x_4)}{f_1(x_4)det(x_1,x_2,x_3)})$$ because the Bloch-Wigner dilogarithm 
satisfies the symmetry $D(z)=-D(\frac{1}{z})$. The same computation shows that also the Bloch-Wigner 
dilogarithms of the second, third, fourth summand in (1) are the negatives of the Bloch-Wigner 
dilogarithms of the second, third, fourth summand in (2). So for each simplex its contribution to the volume from \cite{bfg} is $-\frac{1}{4}$ times its contribution to the volume from \cite{gtz}. Summing over the simplices of the triangulation one obtains the claimed equality.

\end{pf}


\begin{thebibliography}{19}
\bibitem{bas}S.\ Baseilhac: 'Chern Simons theory in dimension three', http://www.math.univ-montp2.fr/$\sim$baseilhac/CS.pdf
\bibitem{bfg} N.\ Bergeron, E.\ Falbel, A.\ Guilloux, 'Tetrahedra of flags, volume and homology of $\mathrm{SL}(3)$', Geom.\ Top.\ 18 (2014), 1911--1971.
\bibitem{bbi}M.\ Bucher, M.\ Burger, A.\ Iozzi: 'The bounded Borel class and complex representations of 3-manifold groups', https://arxiv.org/abs/1412.3428
\bibitem{cs}J.\ Cheeger, J.\ Simons: 'Differential characters and geometric invariants', Geometry and topology (College Park, Md., 1983/84), 50-80, Lecture Notes in Math., 1167, Springer, Berlin, 1985.
\bibitem{fw}E.\ Falbel, Q.\ Wang: 'Duality and invariants of representations of fundamental groups of 3-manifolds into $PGL(3,\C)$.' J.\ London Math.\ Soc.\ (2) 95 (2017), 1-22.
\bibitem{gtz}S.\ Garoufalidis, D.\ Thurston, C.\ Zickert: 'Complex volume of $SL(n,\C)$-representations of 3-manifolds', Duke Math.\ J.\ 164 (2015), 2099-2160.
\bibitem{lod}J.\ Lodder: 'Rigidity of secondary characteristic classes.' Differential Geom.\ Appl.\ 12 (2000), 207-218.
\end{thebibliography}
\end{document}